\documentclass[11pt]{amsart}
\usepackage[T1]{fontenc}
\usepackage[utf8]{inputenc}
\usepackage{amsmath,amssymb,amsthm,mathtools}
\usepackage[a4paper,total={6.5in,8in}]{geometry}
\usepackage{lmodern,microtype}
\usepackage{booktabs}
\usepackage[dvipsnames,svgnames,table]{xcolor}
\usepackage{tikz}
\usetikzlibrary{arrows.meta,calc}
\usepackage[colorlinks=true,linkcolor=RoyalBlue,urlcolor=RoyalBlue,citecolor=PineGreen]{hyperref}
\newtheorem{theorem}{Theorem}
\newtheorem{lemma}{Lemma}
\newtheorem{proposition}{Proposition}

\theoremstyle{remark}

\newcommand{\diam}{\operatorname{diam}}
\newcommand{\conv}{\operatorname{conv}}
\newcommand{\dist}{\operatorname{dist}}
\newcommand{\R}{\mathbb R}
\newcommand{\B}{\mathbb B}
\newcommand{\POLY}{\operatorname{POLY}}
\title[Hausdorff approximation by convex polygons]{Improved perimeter bounds for Hausdorff approximation by convex polygons}
\author{Krishnendu Bhowmick}
\address{Krishnendu Bhowmick \\ TU Chemnitz, Germany.}
\email{krish.combi@gmail.com}
\author{Kurt Klement Gottwald}
\address{Kurt Klement Gottwald \\ TU Chemnitz, Germany.}
\email{kurt-klement.gottwald@math.tu-chemnitz.de}
\date{}
\begin{document}
\begin{abstract}
We investigate how closely a convex polygon with at most $n$ vertices can approximate a plane convex body of perimeter $L$. For free approximation in the Euclidean Hausdorff metric we improve Popov's bound for every $n\ge3$. We also improve his inscribed bound for every $n\ge3$ and his circumscribed bound for $n=3$. The common step is to join successive boundary contacts of the body and estimate the perimeter lost in passing to their convex hull.

\end{abstract}
\maketitle

\section{Introduction}
Fix an integer $n\ge3$. Among plane convex bodies of perimeter one, how large can the error of a best approximating $n$-gon be? The answer depends on whether the polygon is freely placed, must lie in the body, or must contain it. Popov studied all three versions. He conjectured that the regular $(n+1)$-gon is extremal, and obtained upper bounds with the correct leading terms as $n\to\infty$ \cite{Popov,Sendov}. The fixed-$n$ questions are recorded as Problems~2, 3 and~4, respectively, in Brass, Moser and Pach~\cite[\S11.5, p.~468]{BMP}. Here we improve Popov's finite-$n$ bounds by using the geometry of the points where the approximation error is attained.

For nonempty compact sets $K,H\subset\R^2$, their Euclidean Hausdorff distance is
\[
 d_H(K,H)=\max\left\{\sup_{x\in K}\dist(x,H),\ \sup_{y\in H}\dist(y,K)\right\}.
\]
For convex sets the support function gives the equivalent formula
\begin{equation}\label{eq:support}
 d_H(K,H)=\max_{u\in S^1}|h_K(u)-h_H(u)|,\qquad
 h_K(u)=\max_{x\in K}\langle u,x\rangle.
\end{equation}
Let $\POLY_n$ consist of the convex hulls of at most $n$ points. For a nonempty compact convex set $K$, write
\[
 \mathcal E_n(K)=\inf_{P\in\POLY_n}d_H(K,P),\qquad
 \mathcal I_n(K)=\inf_{\substack{P\in\POLY_n\\P\subset K}}d_H(K,P),\qquad
 \mathcal O_n(K)=\inf_{\substack{P\in\POLY_n\\K\subset P}}d_H(K,P).
\]
We call these the free, inscribed and circumscribed errors. A convex body has nonempty interior. We write $L(K)$ for perimeter, counting a segment twice and assigning a point perimeter zero. In the inscribed problem the vertices of $P$ may be chosen anywhere in $K$. When $K$ is a polygon, requiring the vertices of $P$ to be vertices of $K$ gives the subpolygon problem studied by Ahn, Bae, Cheong and Gudmundsson~\cite[Theorem~3]{ABCG}.

Popov's perimeter bounds are
\begin{equation}\label{eq:old}
 \mathcal E_n(K)\le\frac{\tan(\pi/(2n))}{2n}L(K),\quad
 \mathcal I_n(K)\le\frac{\sin(\pi/n)}{2n}L(K),\quad
 \mathcal O_n(K)\le\frac{\tan(\pi/n)}{2n}L(K).
\end{equation}
They are equations~(4), (3) and~(5), respectively, in Popov's Theorem~4~\cite[pp.~69--70]{Popov}. See also Sendov~\cite[\S3.4, Theorem~19]{Sendov} and, for the free bound, Lindenbaum and Bruckstein~\cite[p.~518, (3)]{LB}. Our first result concerns free approximation.

\begin{theorem}\label{thm:main}
For every nonempty compact convex set $K\subset\R^2$ and integer $n\ge3$,
\begin{equation}\label{eq:main}
 \mathcal E_n(K)\le\frac{L(K)}{C_n},
 \qquad
 C_n=n\cot\frac{\pi}{2n}
      +\sqrt{n^2\cot^2\frac{\pi}{2n}+4}.
\end{equation}
The infimum is attained in $\POLY_n$. If $K$ has nonempty interior and does not belong to $\POLY_n$, an optimal polygon can be chosen with exactly $n$ vertices.
\end{theorem}

Since $C_n>2n\cot(\pi/(2n))$, the coefficient in~\eqref{eq:main} is strictly smaller than that in~\eqref{eq:old}. For triangles we obtain a further improvement.
\begin{theorem}\label{thm:triangle}
Every plane convex body $K$ admits a freely placed triangle $T$ such that
\[
 d_H(K,T)\le\frac{L(K)}{1+5\cot(3\pi/20)}.
\]
\end{theorem}
Brass and Lassak~\cite[\S III]{BL} state the conjectured free triangle bound as $L(K)/16$. Table~\ref{tab:triangles} gives the numerical values of the upper bounds proved here.

\begin{theorem}[Inscribed approximation]\label{thm:inner}
For every nonempty compact convex set $K\subset\R^2$,
\[
 \mathcal I_3(K)\le\frac{L(K)}{5\sqrt5-4}.
\]
For every integer $n\ge4$,
\[
 \mathcal I_n(K)\le
 \frac{L(K)}{n\csc(\pi/n)+\sqrt{n^2\csc^2(\pi/n)+4}}.
\]
The infimum defining $\mathcal I_n(K)$ is attained.
\end{theorem}

\begin{theorem}[Circumscribed triangles]\label{thm:outer}
Every plane convex body $K$ admits a circumscribed triangle $T$ satisfying
\[
 d_H(K,T)\le\frac{\sqrt5}{8}L(K).
\]
\end{theorem}

For inscribed triangles the strict improvement follows from the comparison in Section~\ref{sec:inner}. For circumscribed triangles it follows by squaring the two coefficients, since $5/64<1/12$.

\begin{table}[ht]
\centering
\caption{Universal upper coefficients multiplying $L(K)$ for triangles.}\label{tab:triangles}
\small
\begin{tabular}{@{}llll@{}}
\toprule
Approximation & Classical coefficient & Present coefficient & Approximate value\\
\midrule
Free & $\sqrt3/18$ & $1/(1+5\cot(3\pi/20))$ & $0.0924808233$\\[3pt]
Inscribed & $\sqrt3/12$ & $1/(5\sqrt5-4)$ & $0.1392691733$\\[3pt]
Circumscribed & $\sqrt3/6$ & $\sqrt5/8$ & $0.2795084972$\\
\bottomrule
\end{tabular}
\end{table}

In the free problem, Kenderov's alternation theorem provides an optimal polygon with two types of contacts. Its vertices are all at distance $\mathcal E_n(K)$ from $K$, and each edge has a parallel support line of $K$ at the same distance~\cite{Kenderov,KK}. Joining the latter contacts in boundary order produces a second polygon $B\subset K$. We use the contact angles to bound $L(B)$ from below and then use the Hausdorff distance from $B$ to $K$ to estimate $L(K)-L(B)$. Lemma~\ref{lem:deficit} proves the perimeter estimate directly, including for degenerate $B$. A related sectional estimate is due to Carozza, Giannetti, Leonetti and Passarelli di Napoli~\cite[Theorem~1.1]{CGLP}.

For inscribed approximation, a boundary partition with equal errors on its arcs supplies the contacts. For circumscribed approximation, free alternation applied to a parallel body gives a containing polygon whose sides support $K$. For triangles, the three contact chords also control the shape of an inscribed triangle. This extra information yields the bound in Theorem~\ref{thm:outer}.

\section{Balanced contacts and their order}
We state explicitly the consequence of the alternation theorem used in the proof.

\begin{proposition}[Balanced optimal polygon]\label{prop:contacts}
Let $K$ be a plane convex body not in $\POLY_n$, and put $e=\mathcal E_n(K)>0$. There is an optimal convex $n$-gon $P$, with distinct vertices $t_0,\ldots,t_{n-1}$ in counterclockwise order, such that, for the outward unit normal $u_i$ of its edge $t_it_{i+1}$,
\begin{align}
 \dist(t_i,K)&=e,\label{eq:vertexcontact}\\
 h_K(u_i)&=\langle u_i,t_i\rangle+e.\label{eq:edgecontact}
\end{align}
Indices are read modulo $n$.
\end{proposition}

\begin{proof}
The vertices of a minimizing sequence lie in a bounded neighborhood of $K$. Repeating vertices gives $n$-tuples with a convergent subsequence. The convex hull of the limiting tuple attains the infimum. Since $K\notin\POLY_n$, Kenderov's Proposition~3.3 gives exactly $n$ vertices. His Theorem~3.1 gives alternation, and his Definition~2.9 and Corollary~3.0 give the side and vertex equations~\cite{Kenderov}. See also Kenderov and Kirov~\cite[Definition~2.1 and Theorem~2.2]{KK}.
\end{proof}

For this optimal polygon, let $a_i$ be the nearest point of $K$ to $t_i$ and set
\[
 n_i=\frac{t_i-a_i}{e}.
\]
Projection onto a closed convex set is unique. Its supporting-line property gives
\begin{equation}\label{eq:ncontact}
 h_K(n_i)=\langle n_i,a_i\rangle
         =\langle n_i,t_i\rangle-e.
\end{equation}
By~\eqref{eq:support},
\[
 h_P(n_i)\le h_K(n_i)+e=\langle n_i,t_i\rangle\le h_P(n_i).
\]
Thus $n_i$ is an outward normal of $P$ at $t_i$, in the cone bounded by $u_{i-1}$ and $u_i$. It cannot equal either boundary normal:~\eqref{eq:ncontact} and~\eqref{eq:edgecontact} would give support values differing by $2e$. Therefore it lies strictly inside that cone.

\begin{figure}
\centering
\begin{tikzpicture}[x=1cm,y=1cm,>=Stealth,font=\small,
  contact/.style={circle,fill=black,inner sep=1.2pt},
  body/.style={RoyalBlue!85!black,line width=.9pt},
  polygon/.style={black,line width=.9pt},
  normal/.style={->,line width=.65pt},
  translated/.style={->,gray,densely dashed,line width=.65pt},
  error/.style={densely dashed,line width=.65pt}]
% Panel (a): the balanced configuration, with both angles drawn in place.
% e=.6, n_i=(.6,.8), n_{i+1}=(-.6,.8), and u_i=(0,1).
\begin{scope}[scale=.93]
  \coordinate (T) at (3,0);
  \coordinate (U) at (-2.4,0);
  \coordinate (A) at (2.64,-.48);
  \coordinate (C) at (-2.04,-.48);
  \coordinate (B) at (.2,.6);
  \coordinate (Q) at (.2,0);
  \path[fill=RoyalBlue!7]
    (-3,-1.5) .. controls (-2.7,-1.1) and (-2.4,-.75) .. (C)
    .. controls (-1.24,.12) and (-.3,.6) .. (B)
    .. controls (.9,.6) and (1.84,.12) .. (A)
    .. controls (3,-.75) and (3.3,-1.15) .. (3.6,-1.6)
    --(3.6,-2.1)--(-3,-2.1)--cycle;
  \draw[body]
    (-3,-1.5) .. controls (-2.7,-1.1) and (-2.4,-.75) .. (C)
    .. controls (-1.24,.12) and (-.3,.6) .. (B)
    .. controls (.9,.6) and (1.84,.12) .. (A)
    .. controls (3,-.75) and (3.3,-1.15) .. (3.6,-1.6);
  \draw[polygon] (-3.05,-1.3)--(U)--(T)--(3.55,-1.28);
  \draw[polygon,densely dashed] (-3.3,-1.8)--(-3.05,-1.3);
  \draw[polygon,densely dashed] (3.55,-1.28)--(3.77,-1.79);
  \draw[gray,densely dotted] (-.8,.6)--(1.25,.6);
  \draw[error] (A)--(T) node[midway,left=2pt] {$e$};
  \draw[error] (C)--(U) node[midway,right=2pt] {$e$};
  \draw[error] (Q)--(B) node[midway,right=3pt] {$e$};
  \draw[gray,line width=.5pt] (.06,0)--(.06,.14)--(.2,.14);
  \draw[normal] (B)--(.2,1.75) node[above] {$u_i$};
  \draw[normal] (T)--(3.78,1.04) node[above right] {$n_i$};
  \draw[normal] (U)--(-3.18,1.04) node[above left] {$n_{i+1}$};
  % The edge normal is translated to each endpoint without rotation.
  \draw[translated] (T)--(3,1.75) node[above,text=black] {$u_i$};
  \draw[translated] (U)--(-2.4,1.75) node[above,text=black] {$u_i$};
  \begin{scope}[shift={(3,0)}]
    \draw[->,line width=.65pt] (53.130102:.58)
      arc[start angle=53.130102,end angle=90,radius=.58];
    \node at (71.565051:.96) {$\alpha_i$};
  \end{scope}
  \begin{scope}[shift={(-2.4,0)}]
    \draw[->,line width=.65pt] (90:.58)
      arc[start angle=90,end angle=126.869898,radius=.58];
    \node at (108.434949:.96) {$\beta_i$};
  \end{scope}
  \node[contact] at (A) {};
  \node[below left=2pt] at (A) {$a_i$};
  \node[contact] at (C) {};
  \node[below right=2pt] at (C) {$a_{i+1}$};
  \node[contact] at (B) {};
  \node[above left=3pt] at (B) {$b_i$};
  \node[contact] at (T) {};
  \node[above left=3pt] at (T) {$t_i$};
  \node[contact] at (U) {};
  \node[above right=3pt] at (U) {$t_{i+1}$};
  \node[contact] at (Q) {};
  \node[below=4pt] at (Q) {$q_i$};
  \node[RoyalBlue!85!black] at (-2.8,-1.93) {$K$};
  \node at (3.9,-1.4) {$P$};
  \node[align=center] at (.2,-2.65) {(a) Balanced contacts and their angles};
\end{scope}
% Panel (b): moving the left vertex into K gives a general approximation.
% No positive equal-error or optimality assertion is made for this polygon.
\begin{scope}[shift={(8.1,0)},scale=.93]
  \coordinate (Vin) at (-1.8,-.75);
  \coordinate (Vout) at (3,0);
  \path[fill=RoyalBlue!7]
    (-3,-1.5) .. controls (-2.7,-1.1) and (-2.4,-.75) .. (-2.04,-.48)
    .. controls (-1.24,.12) and (-.3,.6) .. (.2,.6)
    .. controls (.9,.6) and (1.84,.12) .. (2.64,-.48)
    .. controls (3,-.75) and (3.3,-1.15) .. (3.6,-1.6)
    --(3.6,-2.1)--(-3,-2.1)--cycle;
  \draw[body]
    (-3,-1.5) .. controls (-2.7,-1.1) and (-2.4,-.75) .. (-2.04,-.48)
    .. controls (-1.24,.12) and (-.3,.6) .. (.2,.6)
    .. controls (.9,.6) and (1.84,.12) .. (2.64,-.48)
    .. controls (3,-.75) and (3.3,-1.15) .. (3.6,-1.6);
  \draw[polygon] (-2.9,-1.76)--(Vin)--(Vout)--(3.55,-1.28);
  \draw[polygon,densely dashed] (-3.1,-1.944)--(-2.9,-1.76);
  \draw[polygon,densely dashed] (3.55,-1.28)--(3.77,-1.79);
  \node[contact] at (Vin) {};
  \node[contact] at (Vout) {};
  \node[RoyalBlue!85!black] at (.2,.23) {$K$};
  \node at (-1.4,-1) {$P$};
  \node[align=center] at (.2,-2.65) {(b) A general, unbalanced approximation};
\end{scope}
\end{tikzpicture}
\caption{Local views of $K$ and $P$.
In (a), $P$ is the balanced optimal polygon. The dashed vertical arrows are translated copies of the edge normal $u_i$, making
$\alpha_i=\angle(n_i,u_i)$ and $\beta_i=\angle(u_i,n_{i+1})$ visible at the two vertices.
Here $q_i$ is the orthogonal projection of $b_i$ onto $t_it_{i+1}$ and
$|b_i-q_i|=|t_i-a_i|=|t_{i+1}-a_{i+1}|=e$.
Panel (b) illustrates a general free approximation with one vertex inside $K$ and one outside. It does not satisfy the balanced-contact conditions of Proposition~\ref{prop:contacts}.}
\label{fig:balanced-contacts}
\end{figure}
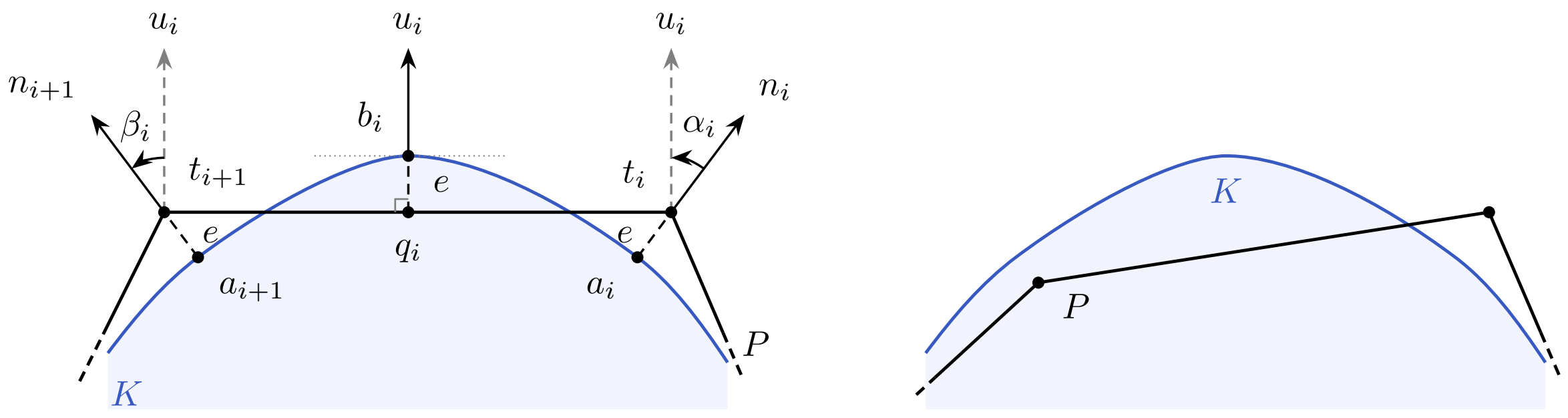

Choose $b_i\in K$ with $\langle u_i,b_i\rangle=h_K(u_i)$. Its projection onto the line of the edge $t_it_{i+1}$ belongs to that edge. Indeed, the distance to the line is $e$, and the distance to $P$ is at most $e$. As the support direction turns counterclockwise, its contact face travels counterclockwise around $\partial K$. Thus the two cyclic orders have the following form.
\[
 n_0,u_0,n_1,u_1,\ldots,n_{n-1},u_{n-1},
 \qquad a_0,b_0,a_1,b_1,\ldots,a_{n-1},b_{n-1}.
\]
At a flat face any choice on the face gives the same weak cyclic order. Consecutive contacts $a_i,b_i$ and $b_i,a_{i+1}$ are distinct by the positive gaps in~\eqref{eq:gap1} and~\eqref{eq:gap2} below. Define the positive normal angles
\[
 \alpha_i=\angle(n_i,u_i),\qquad
 \beta_i=\angle(u_i,n_{i+1}).
\]
Figure~\ref{fig:balanced-contacts}(a) shows the two angles at a single edge.

The normal cones at the vertices give
\begin{equation}\label{eq:angles}
 0<\alpha_i,\beta_i<\pi,
 \qquad \sum_{i=0}^{n-1}(\alpha_i+\beta_i)=2\pi.
\end{equation}
The edge support equations imply the exact gaps
\begin{align}
 \langle u_i,b_i-a_i\rangle&=e(1+\cos\alpha_i),\label{eq:gap1}\\
 \langle u_i,b_i-a_{i+1}\rangle&=e(1+\cos\beta_i).\label{eq:gap2}
\end{align}

\section{Boundary arcs and perimeter deficit}
Define a convex function on $(0,\pi)$ by
\begin{equation}\label{eq:f}
 f(x)=\begin{cases}
       \cot(x/2),&0<x\le\pi/2,\\
       1+\cos x,&\pi/2\le x<\pi.
      \end{cases}
\end{equation}
The values and first derivatives match at $\pi/2$. The second derivatives of the two pieces are $\tfrac12\csc^2(x/2)\cot(x/2)>0$ and $-\cos x\ge0$, respectively.

\begin{lemma}[Boundary arcs]\label{lem:arc}
For the contacts above,
\begin{equation}\label{eq:arc}
 L(K)\ge e\sum_{i=0}^{n-1}\bigl(f(\alpha_i)+f(\beta_i)\bigr).
\end{equation}
\end{lemma}
\begin{proof}
Denote by $l_{a_i}$ the length of the oriented boundary arc from $a_i$ to $b_i$ and by $l_{b_i}$ the length of the oriented boundary arc from $b_i$ to $a_{i+1}$.
Then $L(K)=\sum_{i=0}^{n-1}(l_{a_i}+l_{b_i})$, where $a_n$ is identified with $a_0$.
At almost every point of the oriented boundary arc from $a_i$ to $b_i$, the corresponding unit tangent has projection onto $u_i$ of length at most
\[
 m(\alpha_i)=\begin{cases}
             \sin\alpha_i,&\alpha_i\le\pi/2,\\
             1,&\alpha_i\ge\pi/2.
            \end{cases}=\max\{\sin\varphi:\,\varphi\in[0,\alpha_i]\}.
\]
Integrating this projection and using~\eqref{eq:gap1} gives
\[
 e(1+\cos\alpha_i)=\langle u_i,b_i-a_i\rangle\leq\int_0^{l_{a_i}}m(\alpha_i)ds=l_{a_i}\cdot m(\alpha_i).
\]
Consequently,
\[
 l_{a_i}\geq\frac{e(1+\cos\alpha_i)}{m(\alpha_i)}=e f(\alpha_i).
\]
Projecting the arc from $b_i$ to $a_{i+1}$ onto $-u_i$, identity~\eqref{eq:gap2} similarly gives $l_{b_i}\geq ef(\beta_i)$.
Convex boundaries are rectifiable, have tangents almost everywhere, and have monotonically turning tangent directions. Thus flat sides and nonsmooth points are permitted.
\end{proof}

\begin{lemma}[Perimeter deficit]\label{lem:deficit}
If $B\subset K$ are nonempty compact convex plane sets, with $l=L(B)$, $L=L(K)$, $D=\diam(B)$ and $d=d_H(K,B)$, then
\begin{equation}\label{eq:diameterdeficit}
 L-l\ge\sqrt{D^2+4d^2}-D.
\end{equation}
In particular,
\begin{equation}\label{eq:deficit}
 L\ge\frac{l+\sqrt{l^2+16d^2}}2.
\end{equation}
\end{lemma}
\begin{proof}
The case $d=0$ is perimeter monotonicity. For $d>0$, choose $p\in K$ with $\dist(p,B)=d$ and let $q$ be its projection onto $B$. Translate $q$ to zero, put $u=p/d$, and choose a perpendicular unit vector $v$. Projection gives $\langle u,x\rangle\le0$ for $x\in B$. Write the projection of $B$ onto the $v$-axis as $[-a,b]$. Since $0\in B$, we have $a,b\ge0$, and $a+b\le D=\diam(B)$.

For $u_\theta=u\cos\theta+v\sin\theta$, $|\theta|\le\pi/2$, using $\langle u,x\rangle\cos\theta\leq0$ for $x\in B$, we have
\[
 h_K(u_\theta)\ge d\cos\theta,
 \qquad
 h_B(u_\theta)\le\begin{cases}
                  b\sin\theta,&\theta\ge0,\\
                  a|\sin\theta|,&\theta\le0.
                 \end{cases}
\]
Write $t_+=\max(t,0)$. For $c\ge0$ and $d>0$,
\[
 \int_0^{\pi/2}(d\cos\theta-c\sin\theta)_+\,d\theta=\sqrt{c^2+d^2}-c.
\]
For $c>0$ one integrates up to $\arctan(d/c)$, and for $c=0$ up to $\pi/2$. Since $h_K\ge h_B$, Cauchy's perimeter formula and the preceding support estimates give
\[
 L-l=\int_{-\pi}^{\pi}(h_K-h_B)(u_\theta)\,d\theta
 \ \ge\ \sqrt{a^2+d^2}-a+\sqrt{b^2+d^2}-b .
\]

Since $g(x)=\sqrt{x^2+d^2}-x$ is convex and decreasing for $x\ge0$,
\[
 L-l\ge2g\left(\frac{a+b}2\right)
       \ge2g\left(\frac D2\right)
       =\sqrt{D^2+4d^2}-D.
\]
Finally $D\le l/2$, by considering the two boundary arcs joining diameter endpoints. This also holds for segments and points under the stated conventions. The last expression decreases with $D$, whence
\[
 L-l\ge\sqrt{l^2/4+4d^2}-l/2. \qedhere
\]
\end{proof}

\section{The two contact-angle cases}
Set
\[
 M_n=1+(2n-1)\cot\frac{3\pi}{4(2n-1)}.
\]
We will verify in the next section that
\begin{equation}\label{eq:compare}
 M_n>C_n\qquad(n\ge3),
\end{equation}
which is what we use in the following proof of Theorem~\ref{thm:main}.

\begin{proof}[Proof of Theorem~\ref{thm:main}]
Attainment of the infimum was established in Proposition~\ref{prop:contacts}. If $K\in\POLY_n$, then $\mathcal E_n(K)=0$ and the bound is immediate. This includes all degenerate $K$. Otherwise take the balanced optimal $n$-gon, with error $e>0$.

Suppose first that all $2n$ angles in~\eqref{eq:angles} are at most $\pi/2$. Put
\[
 B=\conv\{b_0,\ldots,b_{n-1}\},\qquad l=L(B).
\]
Let $\tau_i$ be the unit vector obtained by rotating $n_i$ counterclockwise by $\pi/2$. The supporting property of $n_i$ gives
\[
 \langle n_i,b_i-a_i\rangle\le0,
 \qquad\langle n_i,b_{i-1}-a_i\rangle\le0.
\]
In the local frame we have
\[
 u_i=n_i\cos\alpha_i+\tau_i\sin\alpha_i,
 \quad
 u_{i-1}=n_i\cos\beta_{i-1}-\tau_i\sin\beta_{i-1}.
\]
Both cosines are nonnegative. The support gaps~\eqref{eq:gap1}--\eqref{eq:gap2} and the half-angle identity therefore give
\begin{align*}
 \langle\tau_i,b_i-a_i\rangle&\ge e\cot(\alpha_i/2),\\
 \langle\tau_i,a_i-b_{i-1}\rangle&\ge e\cot(\beta_{i-1}/2).
\end{align*}
Therefore
\begin{equation}\label{eq:contactpolygonsingle}
 |b_i-b_{i-1}|\ge e\bigl(\cot(\alpha_i/2)+\cot(\beta_{i-1}/2)\bigr).
\end{equation}
The points $b_i$ occur in cyclic order on $\partial K$, so their consecutive chord lengths sum to the perimeter of $B$. Summing and applying Jensen's inequality gives
\begin{equation}\label{eq:contactpolygon}
 l\ge e\sum_{i=0}^{n-1}\bigl(\cot(\alpha_i/2)+\cot(\beta_i/2)\bigr)
   \ge2n\cot\frac{\pi}{2n}\,e,
\end{equation}
since the angles $\alpha_i$ and $\beta_i$ sum up to $2\pi$. Because $B\in\POLY_n$, global optimality of $e$ implies
\[
 d=d_H(K,B)\ge\mathcal E_n(K)=e.
\]
Lemma~\ref{lem:deficit}, whose lower bound is increasing in both $l$ and $d$, now yields
\[
 L(K)\ge\frac{l+\sqrt{l^2+16d^2}}2\ge C_ne.
\]

In the other case, one contact angle has value $A\ge\pi/2$. Jensen's inequality for the convex function $f$ and Lemma~\ref{lem:arc} yield
\[
 \frac{L(K)}e
 \ge f(A)+(2n-1)f\left(\frac{2\pi-A}{2n-1}\right).
\]
Lastly, the preceding expression is nondecreasing in $A\ge\pi/2$: its derivative is
\[
 f'(A)-f'\left(\frac{2\pi-A}{2n-1}\right)\ge0,
\]
since $f'$ is nondecreasing and $A\ge(2\pi-A)/(2n-1)$. Evaluate at the smallest possible value $A=\pi/2$, where the mean of the remaining angles is $3\pi/(2(2n-1))<\pi/2$. This gives
\[
 f(\frac{\pi}{2})+(2n-1)f\left(\frac{2\pi-\frac{\pi}{2}}{2n-1}\right)=1+(2n-1)\cot\left(\frac{3\pi}{4(2n-1)}\right)=M_n.
\]
Inequality~\eqref{eq:compare} completes this case and the theorem.
\end{proof}

\section{Comparison of the constants}
\begin{lemma}\label{lem:comparison}
For every integer $n\ge3$, $M_n>C_n$.
\end{lemma}
\begin{proof}
We use the elementary inequalities
\begin{equation}\label{eq:cotbounds}
 \cot x<\frac1x\quad(0<x<\pi/2),
 \qquad
 \cot x\ge\frac1x-\frac x2\quad(0<x\le1/2).
\end{equation}
The first follows from $\tan x>x$. For the second, use $\sin x\le x$ and $\cos x\ge1-x^2/2>0$.

First let $n\ge4$, and write $a=n\cot(\pi/(2n))$. Since $a<2n^2/\pi$ and $x\mapsto x+\sqrt{x^2+4}$ increases,
\[
 C_n=a+\sqrt{a^2+4}<\frac{2n^2}\pi+\sqrt{\frac{4n^4}{\pi^2}+4}
     <\frac{4n^2}\pi+\frac\pi{n^2}.
\]
The last inequality follows from \(\sqrt{x^2+4}<x+2/x\) for \(x>0\), applied with \(x=2n^2/\pi\).
The angle $y=3\pi/(4(2n-1))$ is at most $3\pi/28<1/2$. By~\eqref{eq:cotbounds},
\[
 M_n=1+(2n-1)\cot y\ge 1+(2n-1)\left(\frac{1}{y}-\frac{y}{2}\right)=1+\frac{4(2n-1)^2}{3\pi}-\frac{3\pi}8.
\]
Consequently
\begin{equation}\label{eq:nlarge}
 M_n-C_n>
 D(n):=1+\frac{4(n^2-4n+1)}{3\pi}-\frac{3\pi}8-\frac\pi{n^2}.
\end{equation}
The function $D(n)$ increases for real $n\ge4$, since
\[
 D'(n)=\frac{8n-16}{3\pi}+\frac{2\pi}{n^3}>0.
\]
At $n=4$, using $\pi<22/7$ gives
\[
 D(4)=1+\frac4{3\pi}-\frac{7\pi}{16}
     >1+\frac{14}{33}-\frac{11}8
     =\frac{13}{264}>0.
\]
This proves the claim for all $n\ge4$.

For $n=3$ the comparison is $C_3=3\sqrt3+\sqrt{31}<54/5<M_3$.
The last inequality follows by squaring in the identity
$M_3=5\sqrt5-4+5\sqrt{5-2\sqrt5}$.

\end{proof}

\section{A sharper estimate for triangles}
We now retain the diameter in~\eqref{eq:diameterdeficit}, instead of replacing it by half the perimeter. Write
\[
 m=6\sqrt3,\qquad a=4\cot\frac{3\pi}{16},\qquad
 b=1+5\cot\frac{3\pi}{20}=M_3.
\]
The two numerical comparisons needed below are
\begin{equation}\label{eq:trianglecomparisons}
 2+2\sqrt2+a>b,
 \qquad a+\sqrt{(m-a)^2+4}>b.
\end{equation}

\begin{proof}[Proof of Theorem~\ref{thm:triangle}]
If $e=\mathcal{E}_3(K)=0$, then $K\in\POLY_3$, and the assertion is immediate. So we may assume that $e>0$. Scale so that $e=1$. If any contact angle is at least $\pi/2$, the proof of Theorem~\ref{thm:main} already gives $L(K)\ge M_3=b$, which, for $e=1$, is equivalent to the assertion to be proved.

Suppose now that all six contact angles are at most $\pi/2$. Let $B=\conv\{b_0,b_1,b_2\}$, let its three consecutive side lengths be $s_i=|b_i-b_{i-1}|$, and put $l=L(B)$. Set
\[
 \gamma_i=\alpha_i+\beta_{i-1}.
\]
This is the exterior angle of the optimal triangle at $t_i$, so $0<\gamma_i<\pi$ and $\sum_{j\neq i}\gamma_j=2\pi-\gamma_i$. The chord estimates from~\eqref{eq:contactpolygonsingle} and~\eqref{eq:contactpolygon} and convexity of cotangent give
\begin{align}
 s_i&\ge2\cot\frac{\gamma_i}4,\label{eq:oneside}\\
 \sum_{j\ne i}s_j&\ge4\cot\frac{2\pi-\gamma_i}8.\label{eq:othersides}
\end{align}
Consequently
\[
 l\ge F(\gamma_i):=2\cot\frac{\gamma_i}4
                         +4\cot\frac{2\pi-\gamma_i}8.
\]
For $0<x\le2\pi/3$,
\[
 F'(x)=\frac12\left(\csc^2\frac{2\pi-x}8-\csc^2\frac x4\right)\le0,
\]
since $x/4\le(2\pi-x)/8$ and both arguments lie in $(0,\pi/2)$ where $\csc^2(x)$ is monotonically decreasing.

If $l\ge b$, perimeter monotonicity proves the desired bound, since $B\subset K$. Suppose $l<b$. If some $\gamma_i\le\pi/2$, the first comparison in~\eqref{eq:trianglecomparisons} would give
\[
 l\ge F(\gamma_i)\ge F(\pi/2)
   =2+2\sqrt2+a>b,
\]
since $F(x)$ is monotonically decreasing on $(0,2\pi/3]$ as shown above.
This contradicts $l<b$, so every $\gamma_i>\pi/2$. Equation~\eqref{eq:othersides} now shows that the sum of the other two sides is greater than $a$, for each $i$. Hence
\[
 D=\diam(B)=\max_i s_i\le l-a.
\]
Also $d=d_H(K,B)\ge\mathcal E_3(K)=1$, and using $n=3$ and $e=1$,~\eqref{eq:contactpolygon} gives $l\ge m$. Applying~\eqref{eq:diameterdeficit} and monotonicity of $\sqrt{x^2+4}-x$,
\[
 L(K)\ge l+\sqrt{D^2+4}-D
       \ge a+\sqrt{(l-a)^2+4}
       \ge a+\sqrt{(m-a)^2+4}>b.
\]
Here $l\ge m>a$ since $a<4\sqrt3<m$, and the final inequality is~\eqref{eq:trianglecomparisons}.
\end{proof}

\section{Inscribed approximation}\label{sec:inner}
The inner proof uses a balanced boundary partition, but does not require that its polygon be optimal. This distinction allows a direct existence argument.

\begin{lemma}[Equal-error boundary partition]\label{lem:partition}
Let $K$ be a plane convex body not in $\POLY_n$, where $n\ge3$. There are distinct points $a_0,\ldots,a_{n-1}$ in cyclic order on $\partial K$ and a number $e>0$ such that each boundary arc from $a_i$ to $a_{i+1}$ has maximum distance $e$ from the chord $[a_i,a_{i+1}]$. Their convex hull $P$ is an $n$-gon with
\[
 d_H(K,P)=e\ge\mathcal I_n(K).
\]
\end{lemma}
\begin{proof}
Parameterize $\partial K$ counterclockwise by arclength $x:[0,L]\to\partial K$, with $x(0)=x(L)$. For $0\le s\le t\le L$, put
\[
 F(s,t)=\max_{s\le r\le t}\dist\bigl(x(r),[x(s),x(t)]\bigr).
\]
The function \(F\) is continuous on the closed parameter triangle \(0\le s\le t\le L\), since the corresponding boundary arcs and endpoint chords vary continuously in Hausdorff distance. Here the parameter intervals are kept distinct even when their endpoints have the same image under \(x\). In particular, although \(x(0)=x(L)\), the interval \([0,L]\) represents a complete traversal of the boundary and is not identified with a zero-length arc.

For $\lambda=(\lambda_1,\ldots,\lambda_n)$ in the standard simplex, set
\[
 t_0=0,\qquad t_i=L\sum_{j=1}^i\lambda_j,\qquad
 e_i(\lambda)=F(t_{i-1},t_i).
\]
The $e_i$ are continuous, and $e_i=0$ whenever $\lambda_i=0$. They cannot all vanish: otherwise the entire boundary would lie in the union of the $n$ endpoint chords, and $K$ would be their endpoints' convex hull, contrary to $K\notin\POLY_n$. Hence
\[
 f_i(\lambda)=\frac{e_i(\lambda)}{\sum_j e_j(\lambda)}
\]
defines a continuous simplex map preserving every face.

There is a point with $f_i=1/n$ for all $i$. To see this directly from Sperner's lemma, triangulate the simplex with mesh tending to zero. Label every triangulation vertex $v$ by an index $i$ with $f_i(v)\ge1/n$. Such an index exists and satisfies $v_i>0$, so the labelling satisfies the boundary condition of Sperner's lemma. Each triangulation has a small simplex with all $n$ labels. Passing to a convergent subsequence and using continuity gives a point where every $f_i\ge1/n$. Their sum is one, so all equal $1/n$.

The resulting errors therefore equal a common $e>0$. All $\lambda_i$ are positive, and the points $a_i=x(t_i)$, indexed cyclically, are distinct. No support segment of $K$ contains two of these points. Indeed, two consecutive selected points along that segment would be consecutive in boundary order and would enclose a straight arc of error zero. In particular, no three selected points are collinear: a line with three distinct boundary points of a convex body must be a support line. Thus all $a_i$ are vertices of a nondegenerate $n$-gon $P$.

For $z$ on the boundary arc from $a_i$ to $a_{i+1}$, cyclic convexity places $z$ in every inward edge halfplane of $P$ except possibly that of $[a_i,a_{i+1}]$. If $z$ is outside $P$, a segment from $z$ to any $y\in P$ crosses the edge $[a_i,a_{i+1}]$. Consequently
\[
 \dist(z,P)=\dist(z,[a_i,a_{i+1}]).
\]
If $z\in P$, then it lies on that chord and the same equality holds. Taking the maximum over the boundary, and using convexity of distance to $P$ to cover the interior of $K$, proves $d_H(K,P)=e$. Since $P\subset K$, this is at least $\mathcal I_n(K)$.
\end{proof}

For the partition in Lemma~\ref{lem:partition}, choose $b_i$ in the open arc from $a_i$ to $a_{i+1}$ whose distance from its chord is $e$. Let $q_i$ be its projection onto $P$, and put $v_i=(b_i-q_i)/e$. Projection and $K\subset P+e\B$ give
\[
 h_P(v_i)=\langle v_i,q_i\rangle,
 \qquad h_K(v_i)=h_P(v_i)+e=\langle v_i,b_i\rangle.
\]
Thus $v_i$ is a support normal of $K$ at $b_i$, and
\begin{equation}\label{eq:innergaps}
 \langle v_i,b_i-a_i\rangle\ge e,
 \qquad \langle v_i,b_i-a_{i+1}\rangle\ge e.
\end{equation}
If $q_i$ is a chord endpoint, the same inequalities hold even when $v_i$ is oblique to the chord.

Choose any support normal $u_i$ of $K$ at $a_i$. Lift the normal directions in their cyclic order and set
\[
 \alpha_i=\angle(u_i,v_i),\qquad
 \beta_i=\angle(v_i,u_{i+1}).
\]
Every angle is positive, since equality of adjacent normals would contradict~\eqref{eq:innergaps}. The sum of all $2n$ angles is $2\pi$. Define the convex nonincreasing function
\[
 g(x)=\begin{cases}\csc x,&0<x\le\pi/2,\\1,&x\ge\pi/2.\end{cases}
\]
Its left derivative at $\pi/2$ is zero, agreeing with the constant branch. On $(0,\pi/2)$ its second derivative is positive. Projecting the boundary arc from $a_i$ to $b_i$ onto $v_i$, as in Lemma~\ref{lem:arc}, gives length at least $e g(\alpha_i)$. The reversed arc from $a_{i+1}$ to $b_i$ gives $e g(\beta_i)$. The required maximum tangent projections are $\sin\alpha_i$ when $\alpha_i\le\pi/2$, and at most one otherwise. Hence
\begin{equation}\label{eq:innerarc}
 L(K)\ge e\sum_i\bigl(g(\alpha_i)+g(\beta_i)\bigr).
\end{equation}

\begin{proof}[Proof of Theorem~\ref{thm:inner}]
The infimum is attained by compactness of vertex tuples in $K^n$. For a body with interior, its value is unchanged if all vertices are required to lie on $\partial K$. A degenerate inner polygon can first be enlarged to a triangle inside $K$, without increasing the error. For a nondegenerate polygon, extending its vertices along rays from an interior point of that polygon to $\partial K$ gives a boundary-vertex polygon containing the original one, again without increasing the error. If $K\in\POLY_n$, including every degenerate $K$, the error is zero. Otherwise put $e_0=\mathcal I_n(K)>0$ and take the equal-error partition above, whose common error satisfies $e\ge e_0$. Write
\[
 A_n=n\csc\frac\pi n+\sqrt{n^2\csc^2\frac\pi n+4},
 \qquad B_n=1+(2n-1)\csc\frac{3\pi}{2(2n-1)}.
\]

If a contact angle has value $A\ge\pi/2$, Jensen's inequality and~\eqref{eq:innerarc} imply
\[
 \frac{L(K)}e\ge 1+(2n-1)g\left(\frac{2\pi-A}{2n-1}\right)
 \ge B_n.
\]
For the last inequality, note that $g$ is nonincreasing and $3\pi/(2(2n-1))<\pi/2$.
Thus
\[
    \mathcal I_n(K)\leq e\leq\frac{L(K)}{B_n}.
\]

Suppose instead that every contact angle is at most $\pi/2$. Let $\tau_i$ be the positively oriented unit tangent perpendicular to $u_i$. Since $u_i$ supports $K$ at $a_i$, and $b_i,b_{i-1}\in K$,
\begin{equation*}
\langle u_i,b_i-a_i\rangle=\langle u_i,b_i\rangle-h_K(u_i)\le 0,
\qquad
\langle u_i,b_{i-1}-a_i\rangle=\langle u_i,b_{i-1}\rangle-h_K(u_i)\le 0 .
\end{equation*}
Writing $v_i=u_i\cos\alpha_i+\tau_i\sin\alpha_i$, the first inequality in~\eqref{eq:innergaps} and
$\cos\alpha_i\ge 0$ give
\begin{equation*}
e\le\langle v_i,b_i-a_i\rangle
=\cos\alpha_i\,\langle u_i,b_i-a_i\rangle+\sin\alpha_i\,\langle\tau_i,b_i-a_i\rangle
\le\sin\alpha_i\,\langle\tau_i,b_i-a_i\rangle ,
\end{equation*}
whence $\langle\tau_i,b_i-a_i\rangle\ge e\csc\alpha_i$. Similarly, we get $\langle\tau_i,a_i-b_{i-1}\rangle\ge e\csc\beta_{i-1}$. Adding the two bounds
and using $|\tau_i|=1$,
\begin{equation*}
|b_i-b_{i-1}|\ \ge\ \langle\tau_i,b_i-b_{i-1}\rangle
\ \ge\ e\bigl(\csc\alpha_i+\csc\beta_{i-1}\bigr).
\end{equation*}
Set $B=\conv\{b_i\}$ and $l=L(B)$. The cyclic chord estimates and Jensen's inequality now give
\[
 l=\sum_i|b_i-b_{i-1}|
 \ge e\sum_i(\csc\alpha_i+\csc\beta_i)
 \ge 2n\csc(\pi/n)e_0.
\]
The polygon $B$ is also inscribed in $K$ and has at most $n$ vertices, so
$d=d_H(K,B)\ge e_0$. Since the right-hand
side of~\eqref{eq:deficit} is increasing in both $l$ and $d$, Lemma~\ref{lem:deficit} gives, with
$a=n\csc(\pi/n)$ and $l\ge 2ae_0$,
\begin{equation*}
L(K)\ \ge\ \frac{l+\sqrt{l^2+16d^2}}{2}
\ \ge\ \frac{2ae_0+\sqrt{4a^2e_0^2+16e_0^2}}{2}
\ =\ e_0\bigl(a+\sqrt{a^2+4}\bigr)\ =\ A_n e_0 .
\end{equation*}

Combining the two cases proves
\begin{equation}\label{eq:innermin}
 \mathcal I_n(K)\le\frac{L(K)}{\min(A_n,B_n)}.
\end{equation}
Both denominators strictly exceed $2n\csc(\pi/n)$. For $A_n$ this is immediate. For $B_n$, apply strict Jensen to one argument $\pi/2$ and $2n-1$ arguments $3\pi/(2(2n-1))$. Their mean is $\pi/n$, and $g$ is strictly convex between those arguments. Thus~\eqref{eq:innermin} strictly improves the classical bound for every finite $n\ge3$. Lemma~\ref{lem:innerconstants} below simplifies its denominator to the values in the theorem.
\end{proof}

\begin{lemma}\label{lem:innerconstants}
We have $\min(A_3,B_3)=5\sqrt5-4$, and $A_n<B_n$ for every integer $n\ge4$.
\end{lemma}
\begin{proof}
Direct evaluation gives $A_3=4+2\sqrt3>B_3=5\sqrt5-4$. 
We further check $A_4<12<B_4$.

For $0<x\le1$, the inequality $\sin x\ge x-x^3/6$ gives
\[
 \csc x\le\frac1{x(1-x^2/6)}
 =\frac1x+\frac{x}{6-x^2}\le\frac1x+\frac x5.
\]
For $n>4$, put $a=n\csc(\pi/n)\ge n^2/\pi$. Using $\sqrt{a^2+4}\le a+2/a$ for $A_n$ and $\sin y\le y$ for $B_n$ yields
\[
 A_n\le\frac{2n^2}\pi+\frac{2\pi}5+\frac{2\pi}{n^2},
 \qquad B_n\ge1+\frac{2(2n-1)^2}{3\pi}.
\]
Consequently, for $n\ge5$,
\[
 B_n-A_n\ge1+\frac{2n^2-8n+2}{3\pi}
                  -\frac{2\pi}5-\frac{2\pi}{n^2}.
\]
The right side increases with real $n\ge5$ and is positive for $n=5$.
\end{proof}

\section{Circumscribed triangles}\label{sec:outer}
We first obtain an outer balance theorem from Proposition~\ref{prop:contacts}.

\begin{proposition}[Balanced outer polygon]\label{prop:outerbalance}
Let $K$ be a plane convex body and let $e=\mathcal O_n(K)>0$. There is an optimal circumscribed $n$-gon $P$ with vertices $t_i$ and outward edge normals $u_i$ such that
\[
 \dist(t_i,K)=e,
 \qquad h_P(u_i)=h_K(u_i).
\]
\end{proposition}
\begin{proof}
The outer infimum is attained. Indeed, a minimizing sequence lies in a fixed compact parallel neighborhood of $K$. Its vertex tuples have a convergent subsequence, and convex hulls, containment and Hausdorff distance pass to the limit. Put $s=e/2$ and $H=K+s\B$. An outer optimum has
\[
 0\le h_P-h_K\le e,
 \qquad d_H(P,H)\le s.
\]
Conversely, if some $R\in\POLY_n$ satisfied $d_H(R,H)=r<s$, then
\[
 h_R\ge h_H-r=h_K+s-r>h_K,
 \qquad h_R-h_K\le s+r<e.
\]
Thus $R$ would contain $K$ and improve its optimal outer error. It follows that $\mathcal E_n(H)=s$.

Apply Proposition~\ref{prop:contacts} to $H$. The resulting exact $n$-gon satisfies $h_P\ge h_H-s=h_K$, so contains $K$, while its edge equations give $h_P(u_i)=h_H(u_i)-s=h_K(u_i)$. Finally,
\[
 \dist(t_i,H)=s,
 \qquad \dist(x,K+s\B)=\max\{\dist(x,K)-s,0\}
\]
imply $\dist(t_i,K)=2s=e$. Its error relative to $K$ is at most $e$ by the support inequalities, and hence equals $e$.
\end{proof}

\begin{lemma}[Outer contact triangle]\label{lem:outercontact}
Suppose $\mathcal O_3(K)=1$. There are three boundary points of $K$ whose convex hull $B$ has cyclic chord lengths $s_1,s_2,s_3>0$ satisfying
\begin{equation}\label{eq:outerquadratic}
 s_1s_2+s_2s_3+s_3s_1\ge4.
\end{equation}
The assertion allows $B$ to be a segment at this stage, with the limiting perimeter convention.
\end{lemma}
\begin{proof}
Take the triangle $P$ from Proposition~\ref{prop:outerbalance}. Let $a_i$ be the projection of $t_i$ onto $K$ and set $n_i=t_i-a_i$. As in~\eqref{eq:ncontact}, $n_i$ is a support normal of $K$ at $a_i$ and lies strictly inside the vertex normal cone of $P$ at $t_i$. Write
\[
 \alpha_i=\angle(n_i,u_i),\qquad
 \beta_{i-1}=\angle(u_{i-1},n_i),\qquad
 \gamma_i=\alpha_i+\beta_{i-1}.
\]
The $\gamma_i$ are the exterior angles of $P$, so $0<\gamma_i<\pi$ and $\sum_i\gamma_i=2\pi$. Each $\alpha_i$ and $\beta_i$ is at most $\pi/2$. For example,
$a_i\in K\subseteq P$ together with $\langle u_i,t_i\rangle=h_P(u_i)$ gives
$\langle u_i,a_i-t_i\rangle=-\langle u_i,n_i\rangle=-\cos\alpha_i\le0$.
The bound for $\beta_{i-1}$ follows by the same
computation with $u_{i-1}$ in place of $u_i$.

Choose $b_i\in K$ on the supporting line with normal $u_i$, in cyclic order. The support equations give
\[
 \langle u_i,b_i-a_i\rangle=\langle u_i,t_i-a_i\rangle=\langle u_i,n_i\rangle=\cos\alpha_i.
\]

Let $\tau_i$ be the unit vector perpendicular to $n_i$ with
$u_i=n_i\cos\alpha_i+\tau_i\sin\alpha_i$. Expanding the preceding identity in the frame
$(n_i,\tau_i)$ splits it into a normal and a tangential part,
\begin{equation*}
\cos\alpha_i=\langle u_i,b_i-a_i\rangle
=\cos\alpha_i\,\langle n_i,b_i-a_i\rangle+\sin\alpha_i\,\langle\tau_i,b_i-a_i\rangle .
\end{equation*}
The first summand is nonpositive because $\cos\alpha_i\ge0$ and $n_i$ supports $K$ at $a_i$.
Therefore $\langle\tau_i,b_i-a_i\rangle\ge\cot\alpha_i$. The same with $u_{i-1}$ gives
$\langle\tau_i,a_i-b_{i-1}\rangle\ge\cot\beta_{i-1}$.
Consequently
\begin{equation}\label{eq:outerchord}
 s_i:=|b_i-b_{i-1}|\ge\langle\tau_i,b_i-b_{i-1}\rangle\geq\cot\alpha_i+\cot\beta_{i-1}
                    \ge2\cot(\gamma_i/2)>0,
\end{equation}
as depicted in Figure~\ref{fig:lemma7}.
The bound by $2\cot(\gamma_i/2)$ follows from convexity of the cotangent on $(0,\pi/2]$.

\begin{figure}
  \centering
  \definecolor{figblue}{RGB}{55,90,190}
  \begin{tikzpicture}[
    line cap=round, line join=round, x=1.3cm, y=1.3cm,
    dot/.style={circle, fill=black, inner sep=0pt, minimum size=1.6pt},
    arr/.style={-{Latex[length=1.5mm,width=1.1mm]}, line width=0.6pt},
    dim/.style={<->, line width=0.45pt},
    aux/.style={black!45, line width=0.4pt, dash pattern=on 0.8pt off 1.5pt},
    declare function={
      al   = 62;              % alpha_i      (degrees)
      be   = 52;              % beta_{i-1}   (degrees)
      ca   = 1/tan(al);       % cot(alpha_i)     (with e = 1)
      cb   = 1/tan(be);       % cot(beta_{i-1})
      cA   = (tan(al))^2/4;   % curvature of the right branch of the arc
      cB   = (tan(be))^2/4;   % curvature of the left  branch
      xbi  = 2*ca;            % abscissa of b_i
      xbp  = -2*cb;           % abscissa of b_{i-1}
    }
  ]
  %---- K : shading and boundary ---------------------------------
  \fill[figblue!7]
      plot[domain=1.28:0, samples=70]  (\x,{-cA*\x*\x}) --
      plot[domain=0:-1.90, samples=70] (\x,{-cB*\x*\x}) --
      (-1.90,-1.66) -- (1.28,-1.66) -- cycle;
  \draw[figblue, line width=0.85pt] plot[domain=0:1.28, samples=70]  (\x,{-cA*\x*\x});
  \draw[figblue, line width=0.85pt] plot[domain=-1.90:0, samples=70] (\x,{-cB*\x*\x});
  %---- P : the two edges through t_i -----------------------------
\draw[line width=0.6pt] (0,1) -- ({1.42*xbi},{1-2*1.42});
\draw[line width=0.6pt] (0,1) -- ({1.42*xbp},{1-2*1.42});
  %---- supporting line of K at a_i -------------------------------
\draw[black!45, line width=0.4pt] (-2.05,0) -- (1.60,0);
  %---- the segment t_i a_i, of length e ---------------------------
  \draw[black!65, line width=0.5pt, dash pattern=on 2.2pt off 1.7pt] (0,0) -- (0,1);
  %---- n_i and the two edge normals at t_i -----------------------
  \draw[arr] (0,1) -- (0,1.58);
  \draw[arr] (0,1) -- ({0.58*sin(al)},{1+0.58*cos(al)});
  \draw[arr] (0,1) -- ({-0.58*sin(be)},{1+0.58*cos(be)});
  %---- angle arcs ------------------------------------------------
  \draw[line width=0.4pt] (0,{1+0.20}) arc[start angle=90, end angle={90-al}, radius=0.20];
  \draw[line width=0.4pt] ({0.31*cos(90+be)},{1+0.31*sin(90+be)})
        arc[start angle={90+be}, end angle=90, radius=0.31];
  %---- dimension lines for the two spans of length e cot ---------
  \draw[aux] ({ca},0) -- ({ca},-0.58);
  \draw[aux] ({-cb},0) -- ({-cb},-0.58);
  \draw[aux] (0,0) -- (0,-0.62);
  \draw[dim] ({-cb},-0.52) -- (0,-0.52);
  \draw[dim] (0,-0.52) -- ({ca},-0.52);
  %---- the gap between b_{i-1} and b_i ---------------------------
  \draw[aux] ({xbp},{-1.02}) -- ({xbp},{-1.30});
  \draw[aux] ({xbi},{-1.02}) -- ({xbi},{-1.30});
  \draw[dim] ({xbp},{-1.28}) -- ({xbi},{-1.28});
  %---- right angle at a_i ---------------------------------------
  \draw[line width=0.4pt] (-0.09,0) -- (-0.09,0.09) -- (0,0.09);
  %---- points ----------------------------------------------------
  \node[dot] at (0,0) {};
  \node[dot] at (0,1) {};
  \node[dot] at ({xbi},-1) {};
  \node[dot] at ({xbp},-1) {};
  %---- labels ----------------------------------------------------
  \node[anchor=south west, font=\footnotesize] at (0,0) {$a_i$};
  \node[anchor=south east, font=\footnotesize] at (-0.1,0.80) {$t_i$};
  \node[anchor=north east, font=\footnotesize] at ({xbi+0.65},{-0.75}) {$b_{i-1}$};
  \node[anchor=north west, font=\footnotesize] at ({xbp-0.40},{-0.75}) {$b_i$};
  \node[anchor=south, font=\footnotesize] at (0.03,1.62) {$n_i$};
  \node[anchor=south west, font=\footnotesize] at (0.60,1.28) {$u_{i-1}$};
  \node[anchor=south east, font=\footnotesize] at (-0.50,1.38) {$u_i$};
  \node[anchor=west, font=\footnotesize] at (0,1.36) {$\beta_{i-1}$};
  \node[anchor=east, font=\footnotesize] at (0.02,1.43) {$\alpha_i$};
  \node[anchor=north, font=\footnotesize] at (0.4,-0.55)  {$\cot\beta_{i-1}$};
  \node[anchor=north, font=\footnotesize] at (-0.5,-0.58) {$\cot\alpha_i$};
  \node[font=\footnotesize, figblue] at (-1.72,-1.52) {$K$};
  \node[font=\footnotesize] at (0.25,0.27) {$P$};
  \node[anchor=north, font=\footnotesize, align=center] at ({(xbi+xbp)/2},{-1.34}){$s_i$};
\end{tikzpicture}
\caption{The two sides through $t_i$ support $K$ at $b_{i-1}$ and $b_i$. The two short horizontal spans measure the lower bounds on their tangential displacements from $a_i$. In this drawing $e=1$ and the contacts are at the same height, so the bottom span is the chord length $s_i$.}

  \label{fig:lemma7}
\end{figure}
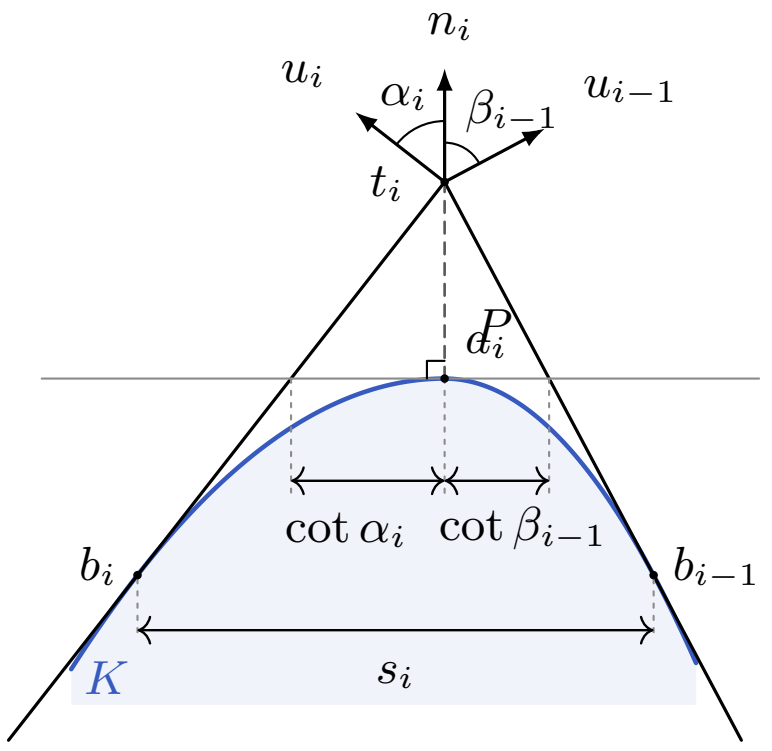

Put $\theta_i=\operatorname{arccot}(s_i/2)\in(0,\pi/2)$. Inequality~\eqref{eq:outerchord} gives $\theta_i\le\gamma_i/2$, and hence $\sum_i\theta_i\le\pi$. The identity
\[
 \sin(\theta_1+\theta_2+\theta_3)
 =\sin\theta_1\sin\theta_2\sin\theta_3
 \left(\cot\theta_1\cot\theta_2+
       \cot\theta_2\cot\theta_3+
       \cot\theta_3\cot\theta_1-1\right)
\]
has a nonnegative left side. As $\cot\theta_i=s_i/2$, the claim \eqref{eq:outerquadratic} follows.
\end{proof}

\begin{lemma}[Parallel outer triangle]\label{lem:paralleltriangle}
Let $B\subset K$ be a nondegenerate triangle with minimum internal angle $A$, and put $d=d_H(K,B)$. Then
\[
 \mathcal O_3(K)\le\frac{d}{\sin(A/2)}.
\]
\end{lemma}
\begin{proof}

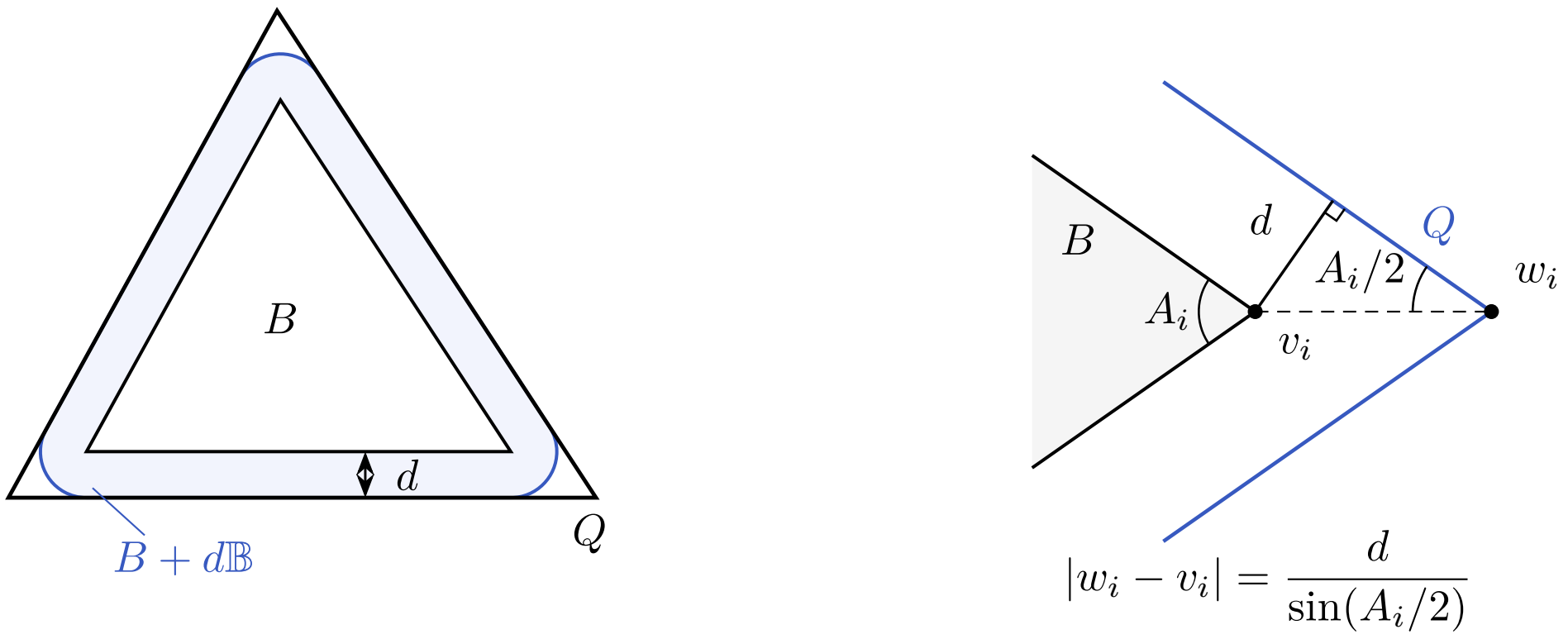
\begin{figure}
\centering
\begin{tikzpicture}[x=1cm,y=1cm,>=Stealth,font=\small,
  dot/.style={circle,fill=black,inner sep=1.2pt}]
\begin{scope}[scale=.97]
  \coordinate (L) at (-1.6,-1);
  \coordinate (R) at (1.9,-1);
  \coordinate (T) at (0,1.9);
  % Exact parallel-body boundary, with d=.38 in drawing coordinates.
  \draw[fill=RoyalBlue!7,draw=RoyalBlue!85!black,line width=.75pt]
    ($(R)+(-90:.38)$)
    arc[start angle=-90,end angle=33.231711,radius=.38]
    --($(T)+(33.231711:.38)$)
    arc[start angle=33.231711,end angle=151.113418,radius=.38]
    --($(L)+(151.113418:.38)$)
    arc[start angle=151.113418,end angle=270,radius=.38]--cycle;
  \draw[fill=white,line width=.75pt] (L)--(R)--(T)--cycle;
  \draw[line width=.9pt] (-2.243654,-1.38)--(2.603260,-1.38)
    --(-.027922,2.636015)--cycle;
  \draw[<->,line width=.55pt] (.7,-1)--(.7,-1.38)
     node[midway,right=3pt] {$d$};
  \node at (0,.1) {$B$};
  \node at (2.55,-1.68) {$Q$};
  \node[RoyalBlue!85!black] at (-.8,-1.88) {$B+d\B$};
  \draw[RoyalBlue!85!black,thin] (-1.12,-1.69)--(-1.55,-1.3);
  \node[align=center] at (.18,-2.8) {(a) Moving each side outward by $d$};
\end{scope}
\begin{scope}[shift={(7.8,.15)},scale=1.05]
  \coordinate (V) at (0,0);
  \coordinate (W) at (1.8,0);
  \coordinate (F) at (.592182,.845723);
  \fill[gray!8] (0,0)--(-1.7,1.19035)--(-1.7,-1.19035)--cycle;
  \draw[line width=.75pt] (-1.7,1.19035)--(V)--(-1.7,-1.19035);
  \draw[RoyalBlue!85!black,line width=.85pt]
    (-.7,1.750519)--(W)--(-.7,-1.750519);
  \draw[densely dashed] (V)--(W);
  \draw[line width=.65pt] (V)--(F) node[midway,above left=1pt] {$d$};
  \draw[line width=.5pt]
    ($(F)+(-35:.11)$)--($(F)+(-35:.11)+(235:.11)$)--($(F)+(235:.11)$);
  \draw[line width=.5pt] (145:.43) arc (145:215:.43);
  \node at (-.67,0) {$A_i$};
  \draw[line width=.5pt] ($(W)+(145:.6)$) arc (145:180:.6);
  \node at ($(W)+(163:1.05)$) {$A_i/2$};
  \node[dot,label={below right:$v_i$}] at (V) {};
  \node[dot,label={above right:$w_i$}] at (W) {};
  \node at (.1,-2.05) {$|w_i-v_i|=\dfrac{d}{\sin(A_i/2)}$};
  \node at (-1.35,.55) {$B$};
  \node[RoyalBlue!85!black] at (1.4,.65) {$Q$};
  \node[align=center] at (0,-2.73) {(b) Displacement at a vertex of angle $A_i$};
\end{scope}
\end{tikzpicture}
\caption{The parallel outer triangle used in Lemma~\ref{lem:paralleltriangle}.
The inclusions $B\subset K\subset B+d\B\subset Q$ make $Q$ an admissible
circumscribed triangle. At a vertex $v_i$ of $B$, the corresponding vertex
$w_i$ of $Q$ lies on the outward angle bisector. The right triangle in (b)
gives the displacement $d/\sin(A_i/2)$.}
\label{fig:parallel-triangle}
\end{figure}

Move each of the three supporting lines of $B$ outward by distance $d$, obtaining a parallel triangle $Q$, depicted in Figure \ref{fig:parallel-triangle} (a). We have $K\subset B+d\B\subset Q$. The vertex of $Q$ corresponding to an angle $A_i$ of $B$ is displaced along the outward angle bisector by $d/\sin(A_i/2)$, as shown in Figure~\ref{fig:parallel-triangle}~(b). This direction lies in the vertex normal cone, so its nearest point in $B$ is the corresponding vertex of $B$. Since distance to a convex set is convex, its maximum over $Q$ occurs at a vertex. Thus
\[
 d_H(Q,K)\le d_H(Q,B)=\frac{d}{\min_i\sin(A_i/2)},
\]
which proves the claim.
\end{proof}

\begin{proof}[Proof of Theorem~\ref{thm:outer}]
If $\mathcal O_3(K)=0$, then $K$ itself is a triangle and the statement is immediate. Otherwise scale so that $\mathcal O_3(K)=1$. Set $C=8/\sqrt5$, and suppose for a contradiction that $L=L(K)<C$. Take $B$ from Lemma~\ref{lem:outercontact}, and relabel its chord lengths as $a\le b\le c$. Write $l=a+b+c=L(B)\le L<C$ and $s=l/2$.

The three distances satisfy the triangle inequalities, even if $B$ is initially degenerate. Thus
\[
 x=s-a,\qquad y=s-b,\qquad z=s-c
\]
are nonnegative and satisfy $x+y+z=s$ and $x\ge y\ge z$. The identity
\[
 ab+bc+ca=s^2+xy+yz+zx,
\]
$x(s-x)\leq\frac{s^2}{4}$ and~\eqref{eq:outerquadratic} give
\begin{equation}\label{eq:outerheron}
 yz\ge4-s^2-x(s-x)\ge4-\frac{5s^2}{4}
                         =4-\frac{5l^2}{16}>0,
\end{equation}
because $l<8/\sqrt{5}$. 
In particular, $B$ is a nondegenerate triangle.

We also have $a>l/6$. Otherwise
\[
 ab+ac+bc\le a(l-a)+\frac{(l-a)^2}{4}
 =\frac{l^2}{4}+\frac{la}{2}-\frac{3a^2}{4}.
\]
The quadratic expression increases for $0\le a\le l/3$, so
\[
ab+ac+bc\le\frac{l^2}{4}+\frac{la}{2}-\frac{3a^2}{4}\le\frac{5l^2}{16}<4,
\]
contradicting~\eqref{eq:outerquadratic}. It follows that $bc\le(l-a)^2/4<25l^2/144$. If $A$ is the minimum angle of $B$, opposite $a$, the half-angle formula and~\eqref{eq:outerheron} therefore yield
\[
 q:=\sin^2(A/2)=\frac{yz}{bc}
 \ge\frac{9(C^2-l^2)}{5l^2}
 \ge\frac{C(C-l)}4.
\]
For the last comparison, divide by $C-l>0$ and use
\[
 \frac{9(C+l)}{5l^2}\ge\frac{18}{5l}
 >\frac{18}{5C}>\frac C4,
\]
where the final inequality is $72>5C^2=64$.

Put $d=d_H(K,B)$. Lemma~\ref{lem:paralleltriangle} and global outer optimality imply $1\le d/\sqrt q$, and hence $d^2\ge C(C-l)/4$. The perimeter deficit inequality~\eqref{eq:deficit} now gives
\[
 L\ge\frac{l+\sqrt{l^2+16d^2}}2
 \ge\frac{l+\sqrt{l^2+4C(C-l)}}2
 =\frac{l+|2C-l|}{2}=C.
\]
This contradicts $L<C$. Rescaling proves the theorem.
\end{proof}

\section{Further questions}
Popov's conjecture asks whether a regular $(n+1)$-gon maximizes the approximation error among bodies of a given perimeter. The free triangle case already exhibits the gap between the conjecture and our bound. For the unit square $Q=[0,1]^2$, the triangle
\[
 T=\conv\left\{\left(\tfrac12,-\tfrac14\right),
                 \left(-\tfrac14,\tfrac34\right),
                 \left(\tfrac54,\tfrac34\right)\right\}
\]
satisfies $d_H(Q,T)=1/4=L(Q)/16$. Indeed, every vertex of either polygon is at distance $1/4$ from the other polygon, and distance to a convex set is convex. The conjecture asks whether \emph{every} plane convex body $K$ admits a triangle at distance at most $L(K)/16$~\cite[\S III]{BL}. Theorem~\ref{thm:triangle} gives $L(K)/(1+5\cot(3\pi/20))$.

For circumscribed polygons, Proposition~\ref{prop:outerbalance} gives support contacts for every $n$. When $n=3$, their three chord lengths control the smallest angle of the contact triangle. This is what allows us to pass from an inscribed polygon to a circumscribed competitor in Lemma~\ref{lem:paralleltriangle}. Can the contacts of a best circumscribed $n$-gon yield an explicit coefficient smaller than $\tan(\pi/n)/(2n)$ for each $n\ge4$?

A separate question concerns the next term as $n$ grows. Let $c_n=\sup_{L(K)=1}\mathcal E_n(K)$ and $p_n=\tan(\pi/(2n))/(2n)$. The regular $(n+1)$-gon gives $p_{n+1}\le c_n$, while~\eqref{eq:old} gives $c_n\le p_n$~\cite[\S3.4, Theorems~19--20]{Sendov}. Since $p_n=\pi/(4n^2)+O(n^{-4})$, these bounds give $n^2c_n\to\pi/4$ and leave a gap of order $n^{-3}$. Theorem~\ref{thm:main} lowers the upper bound for each $n$, but the order of this gap persists. Does $n^3(c_n-\pi/(4n^2))$ have a limit, and if so, what is its value?

\section*{AI involvement in this research}
We used ChatGPT 6 Astra and OpenAI Codex to explore ideas, refine bounds and draft parts of the text. We used Open Math OS to verify parts of the arguments.
All mathematical claims were subsequently verified by the authors, who take full responsibility for the content of this paper.

\end{document}